\documentclass[a4paper,12pt]{amsart} \usepackage{euscript,eufrak,verbatim}
\usepackage[latin1]{inputenc}      
\usepackage[final]{epsfig}
\usepackage{psfrag}

\def\a{\alpha} 
\def\t{\theta}
\def\St{S_{\t}}
\def\e{\varepsilon} 
\def\Z{\mathbb{Z}} 
\def\R{\mathbb{R}} 
\def\S{\mathbb{S}}

\newtheorem{theorem}{Theorem}
\newtheorem{lemma}[theorem]{Lemma}
\newtheorem{proposition}[theorem]{Proposition}
\newtheorem{corollary}[theorem]{Corollary}

\DeclareMathSymbol{\varnothing}{\mathord}{AMSb}{"3F} \begin{document}
\renewenvironment{proof}{\noindent {\bf Proof.}}{ \hfill\qed\\ }
\newenvironment{proofof}[1]{\noindent {\bf Proof of #1.}}{ \hfill\qed\\ }

\title{Periodic billiard orbits\\ in right triangles II}  
\author{Serge Troubetzkoy}
\address{Centre de physique th\'eorique\\
Federation de Recherches des Unites de Mathematique de Marseille\\
Institut de math\'ematiques de Luminy and\\ 
Universit\'e de la M\'editerran\'ee\\ 
Luminy, Case 907, F-13288 Marseille Cedex 9, France}
\email{troubetz@iml.univ-mrs.fr}
\urladdr{http://iml.univ-mrs.fr/{\lower.7ex\hbox{\~{}}}troubetz/} \date{}
\begin{abstract}  
Periodic billiard orbits are dense in the phase space of an
irrational right triangle.  A stronger pointwise density result
is also proven. 
\end{abstract} 
\maketitle
\pagestyle{myheadings}
\markboth{Periodic billiard orbits}{SERGE TROUBETZKOY}
\section{Introduction}
A billiard ball, i.e.~a point mass, moves inside a polygon $P  \subset \R^2$ 
with unit speed along a straight line until it reaches the boundary
$\partial P$, then instantaneously changes direction according to the mirror
law:  ``the angle of incidence is equal to the angle of reflection,'' and
continues along the new line. If the trajectory hits a corner of the polygon,
in general it does not have a unique continuation and thus by definition
it stops there.

It is an open question if there exists a periodic billiard orbit in
every polygon.  None the less for certain classes of polygons one can exhibit
the existence of many periodic orbits. In particular one can ask how dense
the periodic orbits are.  Most of the known results are
about rational polygons, i.e.~polygons for which the angles between the sides
are rational multiples of $\pi$.  The first result in this
direction was that of H.~Masur who showed that the directions of periodic orbits
are dense in a rational polygon \cite{Ma}.  This was strengthened by Boshernitzan
et al.~who showed that for a rational polygon periodic orbits are dense in the 
phase space  \cite{BoGaKrTr}.  In this article they also showed a pointwise density
result in the configuration space:  in a rational
polygon $P$ there exists a dense $G_{\delta}$ set  $G \subset P$ such that for each
point $p \in G$ the orbit of $(p,\theta)$ is periodic for a dense 
subset of directions $\theta \in \S^1$.  Vorobets strengthened 
this result to show
that the set $G$ is also of full measure \cite{Vo}.

Recently I showed that there is an open set $\mathcal{O}$ of right triangles
such that for each irrational $P \in \mathcal{O}$ the set of periodic 
billiard orbits is dense in the phase space \cite{Tr}.  
The main result of this article
is a twofold strengthening of this result.  First of all I extend the result
to all irrational right triangles.

\begin{theorem}\label{phasedense}
Periodic orbits are dense in the phase space of any irrational
right triangle.
\end{theorem}
Remark that this density result also holds for rational
right triangles \cite{BoGaKrTr}. 
Next I strengthen this
density to a pointwise density statement.
\begin{theorem}\label{thm1}
Suppose that $P$ is any irrational right triangle.
Then there exists an at most countable
set $B \subset P$ such that for every
$p \in P\backslash B$ the orbit of $(p,\theta)$ is periodic for a dense 
subset of directions $\theta \in \S^1$.
\end{theorem}

Billiards in right triangles are well known to be equivalent to the
motion of two elastic point masses on a segment (see for example
\cite{MaTa}).  Theorem \ref{thm1} tells us that except for an at
most countable set $B$ of initial positions $0 \le x_1 \le x_2 \le 1$
if $(x_1,x_2)  \not\in B$ then the orbit of $(x_1,v_1),(x_2,v_2)$ is 
periodic for a dense set of velocities $(v_1,v_2)$.

Somewhat surprisingly Theorem \ref{thm1} is stronger than the result of
Vorobets for rational polygons.  There is a special class of rational
polygons known as Veech polygons which are well studied, 
see for example \cite{MaTa} for the definition.
Combing known results on Veech polygons and 
the arguments of this article yields

\begin{proposition}\label{veech}
In $P$ is a Veech polygon then there exists an at most countable
set $B \subset P$ such that for every
$p \in P \backslash B$ the orbit of $(p,\theta)$ is periodic for a dense 
subset of directions $\theta \in \S^1$.
\end{proposition}

Arithmetic or square tiled polygons 
form a subclass of Veech polygons \cite{Zo}.
Let $V(P)$ be the set of corners of $P$. By definition there
are no periodic orbits passing through $p \in V(P)$.
Thus  Theorem \ref{thm1} and Proposition \ref{veech} imply that $B$ is not
empty since $V(P) \subset B$. A simple geometric argument shows

\begin{proposition}\label{arithmetic}
If $P$ is arithmetic then for every
$p \in P\backslash V(P)$ the orbit of $(p,\theta)$ is periodic for a dense 
subset of directions $\theta \in \S^1$.
\end{proposition}
Note that there is
a unique continuation of the billiard orbit through vertices
with angle $\pi/n$. Such a vertex is called regularisable. 
If we consider such orbits as defined by this
continuation then Proposition \ref{arithmetic} holds for $V(P)$ defined
as the nonregularisable vertices.

\section{Strategy}

\begin{figure}
\psfrag*{R0}{{\footnotesize{$R_0$}}}
\psfrag*{R1}{{\footnotesize{$R_1$}}}
\psfrag*{R-1}{{\footnotesize{$R_{-1}$}}}
\psfrag*{1}{{\footnotesize{$1$}}}
\psfrag*{2}{{\footnotesize{$2$}}}
\centerline{\psfig{file=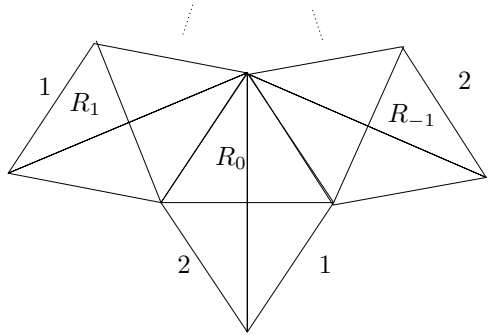,height=44mm}}
\caption{The invariant surface of an irrational triangle.}
\label{fig1}
\end{figure}

Instead of reflection a billiard trajectory in a side of $P$ one can
reflect $P$ in this side and unfold the trajectory to a straight line.
Using this build an invariant flat surface as follows.
Unfold $P$ to a rhombus $R$. Nest unfold $R$ to a
surface $S$ consist of  a
countable union $\{R_n: n \in \Z\}$ of copies of $R$.
The interiors of the $R_n$ are disjoint, and parallel copies of
the same side
are identified by translation (see Figure \ref{fig1} where parallel
sides are labelled by an integer).  
Let $\alpha$ be the smaller angle of the triangle $P$. 
Consider orbits which start in a direction $\theta$ in $R_0$. Then 
the label $n \in \Z$ corresponds
to the billiard orbits in $R$ in the direction $\t_n := \t + 2n \alpha$.
When we fix the initial direction $\theta$ we denote the surface by $\St$.
Sometimes we will refer to the rhombus $R_n$ simply as {\em level $n$}.

Consider the nonsingular 
orbit of a point $x := (s,\theta)$ and the  
sequence of labels $\{l(i): i \in \Z\}$ corresponding to the
sequence of labelled rhombi its orbit visits.
The orbit is called {\it a forward escape orbit} if $\lim_{i \to \infty}
l(i) = \infty$ and a {\it backward escape orbit}  if $\lim_{i \to \infty}
l(i) = -\infty$. Extend these notions to singular escape orbits, 
i.e.~``orbits'' which pass through one (or more) singularity which
are defined by continuity from the left or right.  
The first step is to extend and simplify a result of \cite{Tr}.
\begin{theorem}\label{unique}
In any irrational right triangle $P$, for each $\theta$ which is not 
parallel to the hypotenuse there is
exactly one (possibly singular) forward escape orbit and exactly one (possibly
singular) backward escape orbit on the surface $\St$.
\end{theorem}

Call a direction {\em perpendicular} if it is perpendicular to one of the
legs of the triangle or the hypotenuse. 
Let  $L_\perp$ be the side of the triangle
in question.  An endpoint of $L_\perp$ is {\em good} if the
perpendicular orbit of this point is twice perpendicular, and 
the second perpendicular hit is at an interior point of $L_\perp$.
If $L_\perp$ is a leg is then call $\theta$
{\em end point good } if the endpoint 
which connects $L_\perp$ to the hypotenuse is good.
If $L_\perp$ is the hypotenuse $\theta$ is {\em end point good}
if both end points of $L_\perp$ are good.

\begin{theorem}\label{thmA}
Suppose $P$ is an irrational right triangle.  Fix an end point good 
perpendicular direction $\theta$. 
Then all orbits on $\St$ except for an at most countable collection
of generalized diagonals and the unique forward
and backward escape orbits are periodic.
\end{theorem}
A special case of Theorem \ref{thmA} was proven in \cite{Tr}, namely
when $L_\perp$ is a leg of the triangle and the smaller
irrational angle of $P$ satisfies $\a \in  (\pi/6,\pi/4)$.

Next verify that 
\begin{lemma}\label{lemma1}
Every irrational right triangle has an end point good perpendicular direction.
\end{lemma}

The following Corollary 
follows immediately by combining Theorem \ref{thmA} with
Lemma \ref{lemma1}.

\begin{corollary}\label{thm}
Suppose that $P$ is any irrational right triangle.
Then, for one at least one of the three perpendicular directions $\t$,
the invariant surface $S_{\theta}$ is foliated by periodic orbits except 
for an at most countable collection of orbits/generalized diagonals $O_i$.
\end{corollary}

\begin{proofof} {Theorem \ref{phasedense}} 
The theorem follows from Corollary \ref{thm} and the fact that
the surface $\St$ is dense in the phase space of an irrational
polygon.
\end{proofof}

Turning to the pointwise result we need the following 
\begin{lemma}\label{lemma2}
Suppose $P$ is an irrational polygon. Suppose that 
there exists a direction $\theta \in \S^1$ such
that the invariant surface $S_{\theta}$ is foliated by periodic orbits except 
for an at most countable collection of orbits/generalized diagonals $O_i$. 
Then there exists an at most countable 
set $B \subset P$ such that for every
$p \in P\backslash B$ the orbit of $(p,\theta)$ is periodic for a dense 
subset of directions $\theta \in \S^1$.  
\end{lemma}

\begin{proofof}{Theorem \ref{thm1}}
The theorem follows immediately by combining Corollary \ref{thm} and
Lemma \ref{lemma2}.
\end{proofof}

\section{Proofs}

An orbit segment which begins and ends at a vertex of the polygon
is called a generalized diagonal. A direction $\t$ is called {\it simple}
if there are no generalized diagonals on $\St$.  

\begin{proofof}{Theorem \ref{unique}}
First consider a simple direction $\t$.  For simple directions the Theorem
has already been proven in \cite{Tr}, however for completeness and
to develop the proper notation I sketch the proof here.
The billiard orbits in $R_n$ enter either $R_{n-1}$ or $R_{n+1}$.
A single orbit separates these orbits into two sets, resp.~$R_n^-$ and 
$R_n^+$ which projectively can be thought of as intervals
(see Figure \ref{fig2}).  Call the orbit of such an interval a
{\em strip}.
All the results are obtained by considering compact regions.
Let  $K_N :=R_0^+ \cup R_N^- \cup \cup_{1 \le n \le N-1} R_n$.
Continue the
backward orbit of the
$N-1$ separatrices in $K_N$ until the first time they hit 
the set $R^+_0 \cup R^-_N$. Projectively the set 
$R^+_0 \cup R^-_N$ consists of two interval.  This pull back procedure
splits it into $N+1$ intervals.  Map each of these intervals until
they return to the set $R_0 \cup R_N$.  Note that this set is larger than
$R^+_0 \cup R^-_N$.   This covers $K_N$ by a union
of $N+1$ strips (with disjoint interiors). 
\begin{figure}

\psfrag*{Rn+1}{{\footnotesize{$R_{n+1}$}}}
\psfrag*{Rn-1}{{\footnotesize{$R_{n-1}$}}}
\psfrag*{Rn+}{\footnotesize{{$R_n^+$}}}
\psfrag*{Rn}{{\footnotesize{$R_n$}}}
\psfrag*{Rn-}{\footnotesize{{$R_n^-$}}}
\centerline{\psfig{file=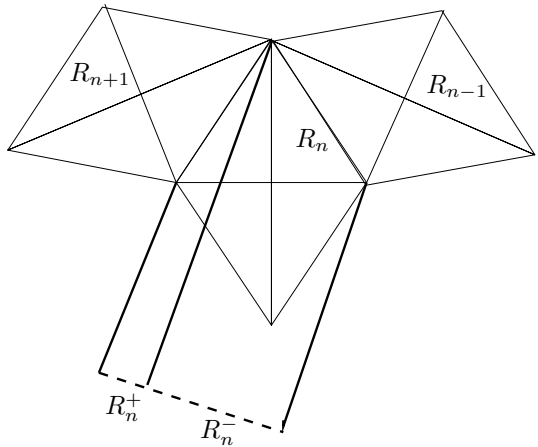,height=58mm}}
\caption{Projective view.}
\label{fig2}
\end{figure}  

Next consider the centers $c_n$ of the rhombi $R_n$.  The orbit $c_n$ has
a centrally symmetric code (labels of the rhombi it visits).  In particular
this implies that each of the $N-1$ centers in $K_N\backslash (R_0^+ \cup R_N^-)$ 
must be in different strips.  By the central symmetry each of the strips which
includes a center must start and stop on the same level 
(i.e.~both on $R_0$ or both on $R_N$)

There are $2 = N+1 - (N-1)$
exceptional strips. A simple argument (Lemma 10 of \cite{Tr}) 
shows that for simple directions there must be an orbit from level $0$ to 
level $N$ and another from $N$ to $0$.  Thus one of the exceptional
strips, call it $E_{0,N}^+$, must start on level 
0 and end on level $N$
while the other, $E_{N,0}^-$, starts on level $N$ and ends on level $0$.
Taking the intersection $\cap_{N \ge 0} E_{N,0}^+$ yields (the forward orbit of)
a unique forward escape orbit. We can replace level $0$ by an arbitrary negative
integer to see that the forward escape orbit does not depend on level $0$,
i.e.~for each $n \in \Z$ 
the image of the forward escape orbit on level $n$ is the forward
escape orbit on level $n+1$.
Similarly $\cap_{N \ge 0} E_{0,N}^-$ yields a unique backward
escape orbit.  

In \cite{Tr} the above  argument
was extended to directions perpendicular to one
of the legs of the triangle formally under the condition that the
angle $\alpha$ of the triangle satisfies $\a \in (\pi/6,\pi/4)$.  
Infact, as mentioned in Section 4 Extensions of \cite{Tr} ``this 
technical assumption guarantees that the orbit starting at the endpoints
of $L$ (the leg of $P$) are simple periodic loops.'' The proof in \cite{Tr}
used an counting argument which showed that even if there are
generalized diagonals the number of  nonexceptional 
strips, i.e.~those which do not contain a point central symmetry, is still
two.

Here I give a much simple proof which holds for all irrational
right triangles for all (nonsimple) directions
except for directions parallel to the hypotenuse.  Note that 
if the direction is
parallel to the hypotenuse the invariant surface breaks into two invariant sets,
the positive levels and the negative levels.  They are 
separated by a generalized
diagonal which is the hypotenuse.  In particular the above picture breaks
down.

Fix a nonsimple direction $\phi$ which is not parallel to a side of the
rhombus.
Suppose that there are at least two distinct 
forward escape orbits on 
$S_{\phi}$.  Let $x_1 := (s_1,\phi)$ and $x_2 := (s_2,\phi)$
be a point on each of these orbits when it passes the last time
through rhombus $R_0$.
Suppose $s_1 < s_2$. Since almost every point is angularly recurrent there
is a nonsingular point $x_3 := (s_3,\phi) \in R_0$ with $s_1 < s_3 < s_2$ whose
orbit returns to $R_0$, i.e.~there is a $M_3 > 0$ such that the angle of
$T^{M_3} x_3$ is $\phi$.  Fix $N > 0$ sufficiently large.
Since $x_1$ and $x_2$ are forward escape orbits
we can find positive integers $M_1,M_2$ such that $T^{M_1} x_1 \in R_{N}$ and 
$T^{M_2}x_2 \in R_{N}$.  Remark that if
the orbit segment $\{x_1,Tx_1,\dots,T^{M_1}x_1\}$, 
is singular, one can
replace $x_1$ by a close point $x_1'$ whose orbit segment has
the same labels at the orbit segment of $x_1$. The same holds for $x_2$. 
In particular the
orbit segment of $x_i'$ starts on level $0$
and arrives at time $M_i$ to level $N$ without revisiting level $0$. This
remark allows us to
assume without loss of generality that the orbits of the $x_i$ are
nonsingular.  
Let $M := \max\{M_1,M_2,M_3\}$.  
Fix $\e > 0$ such that $T^M$ is continuous when restricted to the 
$\e$-ball
around $x_i$, $i=1,2,3$.  By this continuity 
there is a simple directions $\t$
close to $\phi$ such that on the surface $\St$ there are at least two
exception strips of the form $E^+_{0,N}$, a contradiction.  Thus
there is a unique forward escape orbit on $S_{\phi}$.
Similarly there is a unique backward escape orbit.
\end{proofof}

\begin{proofof}{Theorem \ref{thmA}}
Simply repeat the proof  
of Theorem 3 of \cite{Tr}. The fact that the (implicit) assumption 
of this proof is verified is exactly the result of Theorem \ref{unique}.
\end{proofof}

\begin{proofof}{Lemma \ref{lemma1}}
Consider the right triangle and its related rhombus so that the
longer leg is horizontal. Suppose that the angle $\a$ 
between the
hypotenuse and this leg satisfies $2n \alpha < \frac{\pi}{2}$ and
$(2n+1) \alpha >  \frac{\pi}{2}$ for some $n \ge 1$.  Then
the vertical orbit $g$ starting at the left endpoint
of this leg hits the leg again perpendicularly at an interior point
(the case $n=1$, $\alpha \in (\frac{\pi}6,\frac{\pi}{4})$ is illustrated in Figure \ref{fig3}a).  
Thus this leg is end point good for triangles for which  
$\alpha \in \cup_{n \ge 1} (\frac{\pi}{4n+2},\frac{\pi}{4n})$.
\begin{figure}
\psfrag*{L}{{$L$}}
\psfrag*{a}{\footnotesize$\alpha$}
\psfrag*{0}{\footnotesize{$0$}}
\psfrag*{1}{\footnotesize{$1$}}
\psfrag*{g}{{$g$}}
\centerline{\psfig{file=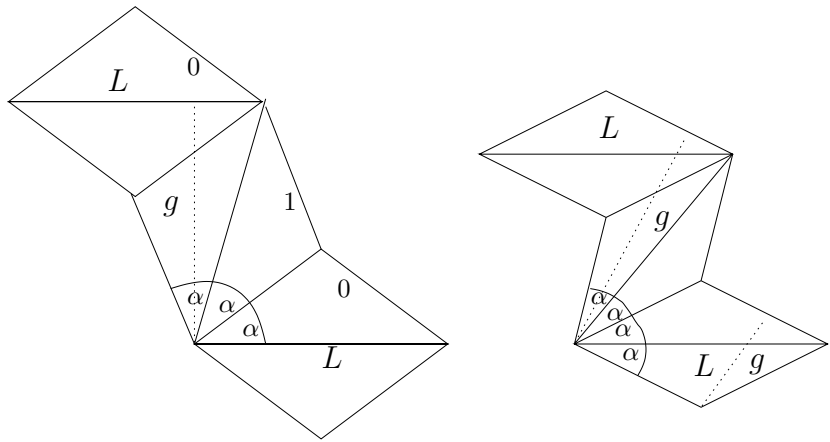,height=58mm}}
\caption{$g$ hits perpendicularly at an interior point.}
\label{fig3}
\end{figure}  

Next consider the orbits perpendicular to the hypotenuse.  If the angle
$\alpha$ satisfies $2n \alpha > \frac{\pi}2$ and $(2n-1) \alpha < \frac{\pi}2$
for some $n \ge 2$ then the perpendicular orbits starting at the endpoints
of the hypotenuse hit the hypotenuse perpendicularly again at an interior
point (the case $n=2$ is illustrated in Figure \ref{fig3}b). Thus the hypotenuse is
end point good for right triangles for which 
$\alpha \in \cup_{n \ge 2} (\frac{\pi}{4n},\frac{\pi}{4n-2})$.
\end{proofof}

\begin{proofof}{Lemma \ref{lemma2}}
Think of $\St$ as being tiled by copies of $P$ which can be enumerated
$P_i$.  Fix $i$ and consider $P_i \cap \{O_j\}$.  This set is an at
most countable collection of parallel line segments $U_{i,j}$.
Let $U_i = \cup_j U_{i,j}$.
Consider the projection $\pi: S_{\theta} \to P$.  This is a countable
to one map. Enumerate the points in $\pi^{-1}(p)$: 
$(p_1,\t_1),(p_2,\t_2),\dots$.
For each $p \in P$ the set of directions $\{\t_k\}$ is dense
in $\S^1$.  By assumption 
each of these directions is the direction 
of a periodic orbit through $p$, unless $p_k$ is in some $U_{i,j}$. 
If the periodic directions through $p$ are not dense then
there exists infinitely many $U_i$ such that $p \in \pi(U_i)$.
Suppose there exists $i,i'$ 
such that $p \in \pi(U_i) \cap \pi(U_{i'})$.
Since the line segments in $\pi(U_i)$ and $\pi(U_{i'})$ are transverse
this intersection is at most countable.  One completes the proof by
taking a union over $i,i'$.
\end{proofof}

\begin{proofof}{Proposition \ref{veech}}
In Veech polygons there is 
a dense set of periodic directions, and each periodic
direction is completely foliated by periodic orbits except for a finite number
of generalized diagonals.  Apply the argument of Lemma \ref{lemma2} to finish
the proof. 
\end{proofof}

\begin{proofof}{Proposition \ref{arithmetic}}
Consider the square, its invariant surface is a torus tiled by
four squares. Fix a rational direction.
Every point in this direction is periodic
except for generalized diagonals.  Suppose that the square unfolds
to a torus.  
Consider the tiling of the plane by  squares representing the torus.  
The intersection points of two generalized diagonals
have rational coordinates.  Thus we conclude that any point in the square
for which at least one coordinate is irrational lies in at most one
generalized diagonal and thus satisfies the conclusions of the proposition.

Consider a rational point $(p_1/q,p_2/q)$ with $gcd(p_1,p_2,q) = 1$.
Tile the plane by $(1/q,1/q)$
squares.  Consider the orbits of slope $a/b$ with
$gcd(a,b)= 1$ starting at the point $(p_1,p_2)$.  Such a
orbit either correspond to a periodic generalized diagonal or a
periodic billiard orbit.  It is a generalized diagonal if and only if
there exists an $i \in \Z$ such that $q|(p_1 + ia)$ and
$q|(p_2 + ib)$ (see Figure \ref{fig4}).  
\begin{figure}
\centerline{\psfig{file=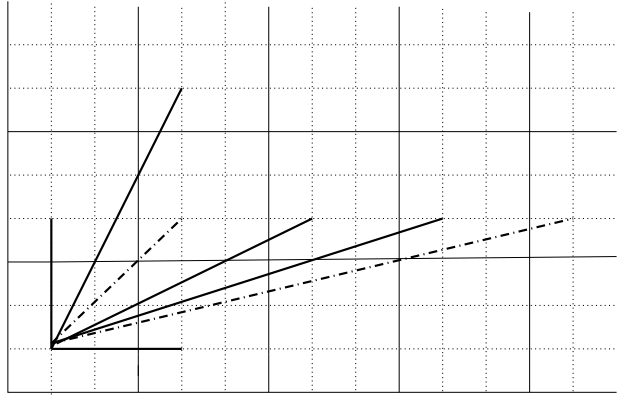,height=53mm}}
\caption{Solid lines are periodic directions while dashed lines are generalized diagonals starting at the point $(1/3,1/3)$.}
\label{fig4}
\end{figure}  
Note that $p_i \in \{0,1,\dots,q\}$ and not being in the set
$V(P)$ is equivalent to 
at least one of the $p_i$
differing from  $0$ or $q$.  If $p_1 \not \in \{0,q\}$ and $q|a$ then
$q$ never divides $p_1 + ia$ and thus
applying the 
above condition yields that the orbit is periodic.
Thus the periodic directions through the given point include the set
$\{qa'/b: qa' \wedge b = 1\}$.
Similarly applying the condition in the case $p_2  \not \in \{0,q\}$ and 
$q| b$ yields that the periodic directions
include the set $\{a/qb': a \wedge qb'= 1\}$.
Both of these sets of directions are dense in $\S^1$.

Now if $P$ is arithmetic then the unfolded surface is a torus cover.
The periodic orbits constructed above for the square lift to periodic
orbits in $P$.
\end{proofof}

\end{document}